\documentclass[10pt]{article}
\usepackage{amssymb}
\font\elevenss=cmss11

\font\eightss=cmss8

\font\sixss=cmss8 at 6pt

\newfam\ssfam
\textfont\ssfam=\elevenss \scriptfont\ssfam=\eightss
 \scriptscriptfont\ssfam=\sixss
\newtheorem {thm}{Theorem}[section]
\newtheorem {lem}[thm]{Lemma}

\def\Cox{\hfill \square}

\def\ee{\epsilon}

\def\E{{\mathbb{E}}}

\def\P{{\mathbb{P}}}

\def\F{{\cal{F}}}

\def\|{\, | \, }
\def\one{{\bf 1}}
\def\Wt{\tilde{W}}
\def\Ut{\tilde{U}}

\begin{document}

\begin{titlepage}
\begin{center}
{\large \bf Network formation by reinforcement learning: 
the long and medium run} \\
\end{center}
\vspace{5ex}
\begin{flushright}
Robin Pemantle \footnote{Research supported in part by National
Science Foundation grant \# DMS 0103635}$^,$\footnote{Department
of Mathematics, The Ohio State University, 231 W. 18 Ave.,
Columbus, OH 43210}
  ~\\
Brian Skyrms \footnote{ School of Social Sciences, University of
California at Irvine, Irvine, CA 92607} \\

\end{flushright}

\vfill

{\bf ABSTRACT:} \hfill \break
We investigate a simple stochastic model of social network formation 
by the process of reinforcement learning with discounting of the past.
In the limit, for any value of the discounting parameter, small,
stable cliques are formed.  However, the time it takes to reach 
the limiting state in which cliques have formed is very sensitive to the
discounting parameter.  Depending on this value, the limiting result
may or may not be a good predictor for realistic observation times.
\vfill


\end{titlepage}

\setcounter{equation}{0}
\section{Introduction} 
\label{ss:intro}

Each day, each member of a small group of individuals selects two others
with whom to interact.  The individuals are of various types, and their 
types determine the payoff to each from the interaction.  That is to say 
that the interaction is modeled as a symmetric 3-person game.  Probabilities 
of selecting individuals evolve by reinforcement learning, where the 
reinforcements are the payoff of the interaction.  We consider two games.  
The first is a degenerate game, ``Three's Company''.  Here there is only 
one type and everyone gets equal reinforcement for every interaction.  
The analysis of  ``Three's company'' is then used in the analysis of 
a second game, a three-person Stag Hunt.  Here there are two types, 
Stag Hunters and Hare Hunters.  Hare Hunters always get a payoff of 3, 
while a Stag Hunter gets a payoff of 4 if he interacts with two other 
Stag Hunters, otherwise he gets nothing.

The point of this modeling exercise is threefold.  First, there is a 
substantial literature, reviewed in the next section, that compares
learning models to laboratory data in order to make inferences about
underlying psychological mechanisms for learning and strategy formation.
Our model explores subgroup formation in a way that is informed 
by previous research on plausible mechanisms and parameter values.  
Secondly, the notion of modeling the co-evolution 
of interaction networks and strategies was set out in~\cite{SP00}.  This
begs the formulation and analysis of basic stochastic network models.
To this end, we provide an analysis of two such models.  Our analysis
is sufficiently robust so as to shed light on similar models, many of 
which have been studied in less depth.  Finally, many of the models
previously studied share a common mathematical description.  We aim 
to provide a collection of rigorous results on this class of models,
which will allow scientists quickly to understand the long- and medium-term
behavior of these stochastic models.

\setcounter{equation}{0}
\section{Reinforcement learning and reinforcement processes}
\label{ss:lit review}

{\em Reinforcement learning}, as the term is used in psychology, 
and  {\em reinforcement models}, as used in applied mathematics, 
are not coextensive.  We will model reinforcement learning by a certain 
kind of reinforcement process, following~\cite{Her70} and~\cite{RE95}.   
Reinforcement learning is to be contrasted to cognitive or belief-based
learning.  In the latter, each agent is assumed to have 
some kind of internal model for the other agents and any unknown
parameters.  Observations alter an agent's beliefs, which then 
govern the agent's future actions.  In a reinforcement learning
model, observations affect future actions through a similar kind
of change of the agent's state, but the translation of observations
into propensities for action are not assumed to be mediated by anything
capable of high-level processing.  Instead, it is assumed that the agent's
actions are governed by previous observations in a simple and explicit way.
One may view this as an exreme form of bounded rationality.
This kind of model is natural when modeling agents with less advanced 
brains than humans have (e.g., slime molds~\cite{OS97}), but has also 
been proposed for humans in situations of low information.  

Stochastic models of reinforcement learning were introduced by 
Estes~\cite{Est50} and by Bush and Mosteller~\cite{BM55}.  
On each of a series of trials the learning agent chooses between a finite
set of alternative acts, and gets a payoff.  The choice is assumed to
be governed by a probability vector that evolves according to 
prescribed learning dynamics.  In the learning dynamics investigated
by Estes and by Bush and Mosteller, the new probability vector after
a trial is a weighted average of the previous probability vector and the
vector putting unit weight on the act just chosen.  The weight on the
unit vector is the product of the payoff and a learning rate parameter.
This kind of learning satisfies what Bush and Mosteller call
{\em independence of path}: the probaiblity vector at time $n+1$ 
depends only on the probability vector at time $n$, the act chosen,
and the magnitude of the payoff.  As a consequence, if the payoffs 
to the acts are fixed, the learning process is a Markov chain with
the probability vector as its state.  Some Markov models for learning 
have been firmly established in the psychology literature,
appearing in survey texts by the early 1970's~\cite{IT69,Nor72}.

Mathematically, a reinforcement process as defined in the literature
on reinforced random walks \hfill \\ \cite{CD97,Dav90,Dav99,Pem88,Pem92}
need not have this Markov property. The current probability can depend
on the history of choices and payoffs, via summary statistics 
or {\em propensities} associated to the possible actions.  The
probability vector is a function of all the propensities (though
in general not a one-to-one function if the process is not Markovian).
The possibility of using such processes to model reinforcement learning 
was introduced by Luce~\cite{Luc59}.  

Luce considered a range of models for the evolution 
of the propensities.  The payoffs for an action taken might modify 
its propensity multiplicatively, additively, or in some combination.  In 
Luce's gamma model, the new propensity, $v'(i)$, for an action $i$, is the
sum of $\gamma_i$ and the product of the old propensity and a factor,
$\beta_i$,
$$v'(i) = \beta_i v(i) + \gamma_i \, .$$
where both $\beta_i$ and $\gamma_i$ are functions of the payoff 
for action $i$.  Luce investigated his models using the linear
response rule,
$$p(i) = \frac{v(i)}{\sum_j v(j)}$$
which simply normalizes the propensities.  However the separation of 
the questions of propensity evolution and response rule opens the 
possibility of other alternatives such as the {\em logistic response 
rule}:
$$p(i) = \frac{\exp (b v(i))}{\sum_j \exp (b v(j))}$$
with $b$ a learning parameter.  This response rule is used in the
learning models of Busemeyer and Townsend~\cite{BT93} and Camerer and 
Ho~\cite{CH99}.  One advantage to the logistic response rule is that
it allows deterrent reinforcement to be modeled, since probabilities
are proportional to exponentials of prepoensities, thus remain positive
when propensities are allowed to become negative.  

The first expermental corroboration of these models, of which we 
are aware, was by Herrnstein~\cite{Her70}.  Thorndike had previously
proposed the ``Law of Effect'', which Herrnstein quantifies as
the ``Matching Law'': the probability of choosing an action is proportional 
to the accumulated rewards.  Let propensities evolve by adding
payoffs, that is, $\beta_i = 1$ and $\gamma_i$ is zero if the action
$i$ was not taken, and is otherwise equal to the payoff.  If we follow
the linear response rule, we obtain Herrnstein's matching law.  Herrnstein
reports data from laboratory experiments with humans as well as with 
animals, from which one may conclude the broad applicability of the model. 

There is a special case whose limiting behavior is well known.  If each 
action is equally reinforced, the process is mathematically equivalent 
to P\'olya's urn process~\cite{EP23}, with each action represented by a
different color of ball initially in the urn.  The process converges
to a random limit, whose support is the whole probability simplex.
In other words, any limiting state of propensities or probabilities 
is possible.  

In 1960, Suppes and Atkinson~\cite{SA60} introduced interactive reinforcement
to model learning behavior in games.  A number of players choose between 
alternatives as before, but the payoffs to each player now depend on the
acts chosen by all players.  Players modify their choice probabilities
by reinforcement learning dynamics of the Bush-Mosteller type.  If joint
actions fix the payoffs and we take the state of the system to be the vector,
indexed by players, of vectors of choice probabilities, then the
dynamics will be Markovian.  

Insofar as multi-agent reinforcement learning has been studied, it
has been largely in the framework of Suppes and Atkinson.  
Macy~\cite{Mac90,Mac91} applies multi-player reinforcement to 
study collective action problems from a bounded rationality viewpoint.  
Borgers and Sarin~\cite{BS97} draw a connection between multi-agent
Bush-Mosteller dynamics and the {\em replicator dynamics} of 
evolutionary game theory~\cite{MS82}, showing that the two coincide 
in a certain limit.  Perhaps the greatest impulse to this direction
of study was the widely cited 1995 paper of Roth and Erev~\cite{RE95}.
They proposed a multi-agent reinforcement model based on Herrnstein's
linear reinforement and response.  Here and in subsequent 
publications~\cite{ER98,BE98}, they show a good fit with a wide range 
of empirical data.  Limiting behavior in the basic model has recently 
been studied by Beggs~\cite{Beg02} and by Ianni~\cite{Ian02}.

In~\cite{SP00}, both basic and discounted versions of Roth-Erev learning 
are applied to social network formation.  Individuals begin 
with prior propensities to interact with each other, and interactions
are modeled as two-person games.  Individuals have given strategies,
and interactions between individuals evolve by reinforcement learning.
The analysis begins with a series of results on ``Making Friends'', a
network formation model in the special case where the game interaction 
is trivial.  Nontrivial strategic interaction is then introduced, 
and it is shown that the co-evolution of network and strategy depends 
on relative rates of evolution as well as on other features of the model.  

The present work is a natural outgrowth of the investigations begun 
in~\cite{SP00}.  In the richer context of multi-agent interactions,
more phenomena arise, namely clique formation and a meta-stable 
state of high network connectivity for an initial epoch whose length
depends dramatically on the discounting parameter.  In
Section~\ref{ss:discussion} we discuss the implications of these 
features for a wide class of models.

\section{Mathematical background}

Our ultimate goal is to understand qualitative phenomena such as
clique formation, or tendency of the interaction frequencies toward
some limiting values.  The mathematical literature on reinforcement
processes contains results in these directions.  It will be instructive
to review these, and to examine the mathematical classification of
such processes, although we will need to go beyond this level of
analysis to explain the behavior of network models such as Three's
Company on timescales we can observe.

Reinforcement processes fall into two main types, trapping and 
non-trapping.  A process is said to be trapping if there are 
proper subsets of actions for each player such that there is
a positive probability that all players always play from this
subset of actions.  For example, if the repetition of any
single vector $({\bf i})$ of actions (action $i_j$ for player $j$) 
is sufficiently self-reinforcing that it might cause action ${\bf i}$
to be perpetuated forever, then the process is trapping.  The specific
dynamics investigated by Bush and Mosteller in 1955 are trapping, as
are most logistic response models.  By contrast, models that give
all times in the past an equal effect on the present, such as
Herrnstein's dynamics and Roth-Erev dynamics, tend not to be trapping.

One of several modifications suggested by Roth and Erev to maximize
agreement of their model with the data is to introduce a discounting
parameter $x \in (0,1)$.  The past is discounted via multiplication
by a factor of $(1-x)$ at each step.  Formally, this is a version of
Luce's gamma model with $\beta_i = 1-x$ for all $i$.  It is known from 
the theory of urn processes that discounting may cause trapping.
For example, it follows from a theorem of H. Rubin reported 
in~\cite{Dav90} that if P\'olya's urn is altered by discounting
the past, there will be a point in time beyond which only one
color is ever chosen.  This holds as well with Roth-Erev type models:  
the discounted Roth-Erev model is trapping, while the undiscounted model 
is not.  In~\cite{SP00}, discounted and nondiscounted versions 
of several games are studied, and equilibria examined for stability.
Again, discounting causes trapping, and we investigate the robustness
of the trapping when the discounting parameter becomes negligible.
In a related paper, Bonacich and Liggett~\cite{BL02} investigate 
Bush-Mosteller dynamics in a two-person interaction representing 
gift giving.  Their model has discounting, and they find a set of
trapping states.  

It is in general an outstanding problem in the theoretical
study of reinforcement models to show that trapping must occur
with probability~1 if it occurs with positive probability.
This was only recently proved, for instance, for the reinforced
random walk on a graph with three vertices, via a complicated
argument~\cite{Lim01}.  Much of the effort that has gone into
the mathematical study of these models has been directed at these
difficult limiting questions.  In the non-trapping case, even though
the choice of action does not fixate, the probaiblities for some
of the actions may tend to zero.  A series of papers in the 1990's
by Benaim and others~\cite{BH95,Ben98,Ben99} establishes some basic tests 
for whether in undiscounted Roth-Erev type models, probabilties will 
tend toward determinstic vectors.  

From the point of view of applications, in a situation where it
can be proven or surmised that trapping occurs, we are mainly 
interested in characterizing the states in which we may become trapped
and in determining how long will it be before the process becomes 
trapped.  Recalling our initial discussion of modeling goals, we are 
particularly interested in results that are robust as parameters
and modeling details vary, or when they are not robust, of understanding 
how these details of the model affect observed qualitative behavior.

\setcounter{equation}{0}
\section{Three's Company: a ternary interaction model} 
\label{ss:three's company}

\subsection{Specification of the model}

The game ``Three's Company'' models collaboration of trios of
agents from a fixed population.  At each time step, each agent
selects two others with whom to form a temporary collusion.  
An agent may be involved in multiple collusions during a single 
time step: one that she initiates, and zero or more initiated by
another agent.  Analogously to the basic game ``Making Friends'',
introduced in~\cite{SP00}, Three's Company has a constant reward
structure: every collaboration results in an identical positive outcome,
so every agent in every temporary collusion increases by an 
identical amount her propensity to choose each of the other two
agents in the trio.  The choice probabilities follow what could be
called {\em mulitlinear response}.  The probability of an agent choosing
to form a trio with two other agents $i$ and $j$ is taken to
be proportional to the product of her propensity for $i$ with her
propensity for $j$.  In addition to providing a model for self-organization
based on a simple matching law type of response mechanism, this model
is meant to provide a basis for the analysis of games such as the 
three person stag hunting game discussed in the next section.  We
now give a more formal mathematical definition of Three's Company,
taken from~\cite{PS03a}.

Fix a positive integer $N \geq 4$, representing the size of the 
population.  For $t \geq 0$ and $1 \leq i , j \leq N$, define 
random variables $W(i,j,t)$ and $U(i,t)$ inductively
on a common probability space $(\Omega , \F , \P)$ as follows.
The $W$ variables are positive numbers, and the $U$ variables are
subsets of the population of cardinality~3.  One may think of the
$U$ variables as random triangles in the complete graph with a 
vertex representing each agent.  The variable $U(i,t)$ is equal to
the trio formed by agent $i$ at time $t$.  The $W$ variables represent 
propensities: $W(i,j,t)$ will be the propensity for player $i$ to choose
player $j$ on the time step $t$.  The initialization is
$W(i,j,0) = 1$ for all $i \neq j$, while $W(i,i,0) = 0)$.  
We write $W(e,t)$ for $W(i,j,t)$ when $e$ is the edge (unordered set)
$\{ i , j \}$ (note that the evolution rules below imply that
$W(i,j,t) = W(j,i,t)$ for all $i, j$ and $t$).  The
inductive step, for $t \geq 0$, defines probabilities (formally, 
conditional probabilities given the past) for the variables $U(i,t)$ 
in terms of the variables $W(r,s,t)$, $r , s \leq N$, and then defines 
$W(i,j,t+1)$ in terms of $W(i,j,t)$ and the variables $U(r,t)$, $r \leq N$.
The equations are:
\begin{eqnarray}
\P (U(i,t) = S \| \F_t) & = & \frac{\one_{i \in S} \prod_{r,s \in S, r < s}
   W(r,s,t)}{\sum_{S' : i \in S'} \prod_{r,s \in S' , r < s}
   W(r,s,t)} \, ; \label{eq:probs} \\[1ex]
W(i,j,t+1) & = & (1-x) W(i,j,t) + \sum_{r=1}^N \one_{i,j \in U(r,t)} \, .
   \label{eq:updates}
\end{eqnarray}
Here $(1-x)$ is the factor per unit time by which the past is discounted,
and the $\sigma$-field conditioned on is the process up to time $t$,
$$\F_t := \sigma \left \{ W(i,j,u) : u \leq t \right \} \, .$$

The following alternative statement of the evolution 
equation~(\ref{eq:probs}) is useful for those familiar with the 
analytic machinery (c.f\ \cite{Pem92}) that is typically used to 
reduce such a process to a stochastic approximation.  Think of the 
normalized matrix 
$${\bf W}_t := \frac{1}{\sum_{i,j} W(i,j,t)} \, W(\cdot , \cdot , t)$$ 
as the state vector.  This is then an asymptotically time-homogeneous 
Markov chain, with an evolution rule
\begin{equation} \label{eq:state}
\E \left ( {\bf W}_{t+1} - {\bf W}_t \| \F_t \right ) = 
   g(t) \left [ \mu ({\bf W}_t) + \xi_t \right ] 
\end{equation}
where $g(t) = 1/x + O(1/t)$, the drift vector field $\mu$ maps the
simplex of normalized matrices into its tangent space and may 
be explicitly computed, and $\xi_t$ are martingale increments of order~1.
In the non-discounted case, $g(t) = 1/t + O(1/t)$, and much information 
about the long term behavior of this process can be discovered by 
an analysis of the the flow $d{\bf X} / dt = \mu ({\bf X})$~\cite{Ben99}.  
In the discounted case, $g(t)$ does not go to zero and an alternative
analysis is required.

\subsection{Analysis of the model}

Equations~(\ref{eq:probs}) and~(\ref{eq:updates}) completely specify
the model for the given parameters $N$ and $x$.  Simulations for a
population of size 6 ($N=6$) showed the following behavior.

When $x = .5$ (a rather steep discount rate, though not unheard of in
psychological laboratory experiments~\cite{BS02}), all 1,000 trials
broke up into two cliques of size 3, with no interactions across
clique boundaries.  In larger populations, with the same discount rate, 
again decomposition into cliques occurs, this time of sizes 3, 4 and 5, 
whose members interact exclusively with other members of the same clique.  

When $N = 6$ and $x = .4$ we found that 994 out of the 1,000 trials
had decomposed into two cliques of three (we allowed the process to
continue for 1,000,000 time steps).  When $x$ was decreased to $.3$,
only 13 of the 1,000 trials showed decomposition into cliques, while
in the remainder of the trials all six members of the population 
remained well connected through the 1,000,000 time steps.  Finally,
when $x = .2$, a reasonable discount rate for individuals though
still steeper than in most economic models, none out of 1,000 trials
had broken into cliques.  All six members of the population remained
well connected after 1,000,000 time steps.

To summarize the simulation data, high discount rates lead to
trapping, with each agent restricting her choices to members of
a clique of size~3 (or, in larger populations, size~4 or~5).  Less
steep discount rates lead to less trapping or no trapping at all.
Interestingly, the simulation data is contradicted by the following
theorem, proved in the appendix.

\begin{thm} \label{th:limit}
In Three's Company, with any population size $N \geq 6$ and any
discount rate $x \in (0,1)$, with probability~1 the population may 
be partitioned into subsets of sizes~3, 4 and~5, such that each member 
of each subset chooses each other with positive limiting frequency, 
and chooses members outside the subset only finitely often.  Every partition 
into sets of sizes~3, 4 and~5 has positive probability of occurring.
\end{thm}

In other words, despite the simulation data, trapping always occurs.
The set of traps is the set of all ways of decomposing into cliques 
of sizes~3,~4 and~5.
The apparent contradiction between the simulation and the theorem 
is resolved by Theorem~\ref{th:transience time}, whose proof
is given in the companion paper~\cite{PS03a}.  The theorem states 
that the time for the population to break into cliques increases
exponentially in $1/x$ as the discount rate $1-x$ increases to~1.    

\begin{thm} \label{th:transience time}
For each $N \geq 6$ there is a $\delta > 0$ and numbers $c_N > 0$ such that 
in Three's Company with $N$ players and discount rate $1-x$,
the probability is at least $\delta$ that each player will play with each
other player beyond time $\exp (c_N x^{-1})$.   $\Cox$
\end{thm}

\section{The three player stag hunt}

\subsection{Specification of the model}

We now replace the uniformly positive reward structure by a nontrivial
game, which is a three player version of Rousseau's Stag Hunt.  For
the purposes of our model, agents are divided into two types, hare
hunters and stag hunters.  That is, we model strategic choice as 
unchanging, at least on the time scale where network evolution is
taking place.  No matter which other two agents a hare hunter goes
hunting with, the hare hunter comes back with a hare (hares can
be caught by individuals).  A stag hunter, on the other hand, comes
home empty-handed unless in a trio of three stag hunters, in which 
case each comes home with one third share of a stag.  One third of a stag
is better than a whole hare, but evidently riskier because it will
not materialize if any member of the hunting party decides to play it safe
and focus attention on bagging a hare.  In the three player stag hunting 
game, as in Three's Company, at each time step each agent chooses two
others with whom to form a collusion.  The payoffs are as follows.
Whenever a hare hunter is a member of a trio, his reward is 3.  A stag
hunter's reward is 4 if in a trio of three stag hunters and 0 otherwise. 
A formal model is as follows.

Let $N = 2n$ be an even integer representing the size of the population
and let $x \in (0,1)$ be the discount parameter.
The variables $\{ W(i,j,t) , U(i,t) : 1 \leq i , j \leq N ; t \geq 0 \}$
are defined again on $(\Omega , \F , \P)$ with the $W$ variables taking
positive values and representing propensities and the $U$ variables
taking values in the subsets of $\{ 1 , \ldots , N \}$ of cardinality~3 
and representing choices of trios.  We initialize the $W$ variables
by $W(i,j,0) = 1 - \delta_{ij}$, just as before, and we invoke a linear
response mechanism~(\ref{eq:probs}) just as before.  Now, instead of
the trivial reward structure~(\ref{eq:updates}), the propensities evolve
according to the hunting bounties
\begin{eqnarray} 
W(i,j,t+1) & = & (1-x) W(i,j,t) 
   + 3 \one_{i \leq n} \sum_{r=1}^N \one_{i,j \in U(r,t)} \nonumber \\[2ex]
&& + 2 \one_{i > n} \sum_{q,r,s=n+1}^N  
   \one_{i \in U(q,t) = \{ q,r,s \} } \, .  \label{eq:bounties}
\end{eqnarray}
The factor in front of the last sum is~2 rather than~4 because 
the sum counts the trio $\{ q,r,s \}$, chosen by agent $q$, 
exactly twice: as $(q,r,s)$ and as $(q,s,r)$.  

\subsection{Analysis of the model}

The propensities for stag hunters to choose rabbit hunters remain
at their initial values, whence stag hunters choose other stag 
hunters with limiting probability~1.  The stag hunters are never
affected by the rabbit hunters' choices, so the stag hunters mimic
Three's Company among themselves precisely except for the times, 
numbering only $O(\log t)$ by time $t$, when they choose rabbit hunters.  
We know therefore, that eventually they fall into cliques of size
3, 4 and 5, but that this will take a long time if the discount
parameter is small.

Rabbit hunters may form cliques of size 3, 4 and 5 as well, but 
because they are rewarded for choosing stag hunters, they may also
attach to stag hunters.  The chosen stag hunters have cliques of their 
own and ignore the rabbit hunters, except during the times that 
they are purposelessly called to hunt with them.  These attachements
can be one rabbit continually calling on a particular pair of stags
or two rabbits continually calling on a single stag.  In either case
the one or two rabbits are isolated from all hunters other than their
chosen stag hunters.  

What matters here is not the details of the trapping state but the
time scale on which the trap forms and the likelihood of a rabbit
hunter ending up in a sub-optimal trap\footnote{In this model, since any
number of collusions is permitted for an agent on each time step,
the sub-optimality is manifested not through wasted time on the
stag hunter's part.  Instead, it is a societal opportunity cost,
borne by all the rabbit hunters passed over in favor of the stag hunter.}.
This likelihood decreases as the discount rate becomes small for the
following reason.  Rabbit hunters choosing to hunt with stag hunters are 
getting no reciprocal invitations, whereas any time they choose
to hunt with other rabbit hunters, their mutual success creates a
likelihood of future reciprocal invitations.  These reciprocal invitations
are then successful and increase the original hunter's propensity for
choosing the other rabbit hunter.  Thus, on average, propensity for
a rabbit hunter to form a hunting party with other rabbit hunters
will increase faster than propensity to call on stag hunters, and 
the relative weights will drift toward the rabbit-rabbit groupings.
The smaller the discount parameter, $x$, the more chance this has to
occur before a chance run of similar choices locks an agent into a particular 
clique.  

Simulations show that stag hunters find each other rapidly. With 6 stag
hunters and six rabbit hunters and a discount rate of .5, the probability
that a stag hunter will visit a rabbit hunter usually drops below half
a percent in 25 iterations. For 50 iterations of the process this always
happened in 1000 trials, and this remains true for values of x between .5
and .1. For x=.01, 100 iterations of the process suffices for stag hunters
to meet stag hunters at this level and for 200 iterations are enough when
x=.001. Rabbit hunters find each other more slowly, except when they are
frozen into interactions with stag hunters. When the past is heavily
discounted the latter possibility a serious one. At x=.5, at least one
rabbit hunter interacted with a stag hunter (after 10,000 iterations) in 384
of 1,000 trials. This number dropped to 217 for x=.4, 74 for x=.3, 6 for x=.2,
and 0 for x=.1. Reliable clique formation among stag hunters is much slower
in line with results of the last section, taking about 100,000 iterations for
x=.5 and 1,000,000 iterations for x=.4.

\subsection{Further discussion}
\label{ss:discussion}

The two models discussed in this paper are highly idealized.  But from 
these, we can learn some general principles as to how to analyze a much 
wider class of models.  

The first principle is that when $x$ is near zero, the process should
for a long time behave similarly to the non-discounted process ($x=0$).
Here, following~\cite{Ben99}, one must find equlibria for the flow
$d {\bf X} / dt = \mu ({\bf X})$, and classify these as to stability.
Unstable equilibria, in general, do not matter (though see~\cite{PS03b}
for cases in which the effects of unstable equilibria may last quite
a while).  Stable equilibria may be possible trapping states, or may
not be.  The interesting case is when a stable equilibrium for the 
non-discounted process is not a possible trapping state for the
discounted process.  In this case, the process may get pseudo-trapped
there, that is, may remain there for a very long time.  Just how long 
will depend on the model, though Theorem~\ref{th:transience time} 
extends rather robustly to a broader class of linearly stable states 
(for the non-discounted process) that are non-trapping for the discounted 
process.

Another mathematical technique relevant to these analyses, which we 
have not yet tried to apply, is quasi-stationary analysis.  Recall that
equation~(\ref{eq:state}) describes and asymptotically time-homogeneous 
Markov chain.  If there is trapping, this chain is not ergodic.  A chain
that is not ergodic may be conditioned to stay in a set of transient
states.  The stationary measure of the conditioned chain is called a
quasi-stationary measure for the original chain.  The study of these
was begun in the 1960's by Seneta and others (see, e.g.,~\cite{DS65})
and there is now an extensive literature.  In particular, it is sometimes
possible to understand the time scale on which the process leaves the
transient states.  

\section{Conclusion}
\label{ss:conclusion}

Our analysis reinforces the emphasis of Suppes and Atkinson, and of 
Roth and Erev, on the medium run for empirical applications.  Long run 
limiting behavior may simply never be seen.  It is useful to quantify
the time scale on which we can expect medium run behavior to persist,
and Theorem~\ref{th:transience time} is meant to serve as a prototypical
result in this direction.  Indeed, Theorem~\ref{th:transience time} is
proved via a stronger result~\cite[Theorem~4.1]{PS03a}, which applies
to many trapping models as the discount rate becomes negligible. 
As to the nature of the medium run behavior, analyses tend to be 
model-dependent.

\setcounter{equation}{0}
\section{Appendix: proof of Theorem~\protect{\ref{th:limit}}}
\label{ss:appendix}

Let $G(t)$ be the graph whose 
edges are all $e$ such that $e \subseteq U(i , t)$ for some $i$, 
that is, the set of edges whose weights are increased from 
time $t$ to time $t+1$.  The following two easy lemmas capture 
some helpful estimates.  

\begin{lem} \label{lem:min}
\begin{enumerate}
\item $$\sum_e W(e,t) \to 3 N x^{-1}$$
exponentially fast as $t \to \infty$.
\item If $e \in G(t)$ then 
$$W(e,t+k) \geq (1-x)^{k-1} \, .$$ 
\end{enumerate}
\end{lem}

\noindent{\sc Proof:}  The first part is a consequence of the equation
for the evolution of the total weight:
$$\sum_e W(i,t+1) = (1-x) \sum_e W(e,t) + 3N \, .$$
The second part follows from the first, and from the fact that when
$e \in G(t)$ then $W(e,t+1) \geq 1$ and hence $W(e,t+k) \geq(1-x)^{k-1}$.
$\Cox$

Let $\overline{G}$ denote the transitive (irreflexive) closure of 
a graph $G$; thus $\overline{G}$ is the smallest disjoint union 
of complete graphs that contains $G$.  

\begin{lem} \label{lem:completion}
There are constant $c$, depending on $N$ and $x$, such that 
with probability at least $c$, every edge $e \in \overline{G(t)}$
satisfies 
\begin{equation} \label{eq:min}
W(e,t+N^2) \geq (1-x)^{N^2} \, .
\end{equation}
\end{lem}

\noindent{\sc Proof:} Let $H$ be a connected component of $G(t)$.
Fix a vertex $v \in H$ and let $w$ be any other vertex of $H$.
There is a path from $v$ to $w$ of length at most $N$; denote
this path $(v = v_1 , v_2 , \ldots , v_r = w)$.  If $r=2$ then  
the inequality~(\ref{eq:min}) for $e \in G(t)$ follows from
Lemma~\ref{lem:min}.  If $r \geq 3$, we let $E(H,v,w,1)$ be the
event that for every $2 \leq j \leq r-1$, the edge $\{ v_{j-1} , v_{j+1} \}$
is in $G(t+1)$.  Since this event contains the intersection over $r$ of
the events that $U(v_j , t) = \{ v_j , v_{j-1} , v_{j+1} \}$, 
since Lemma~\ref{lem:min} bounds each of these probabilities from
below, and since the events are conditionally independent given $\F_t$,
we have a lower bound on the probability of $E(H,v,w,1)$.  In general,
for $1 \leq k \leq r-2$, let $E(H,v,w,k)$ be the event that for every 
$2 \leq j \leq r-k$, the edge $\{ v_{j-1} , v_{j+k}$ is in $G(t+k)$.
We claim that conditional on $E(H,v,w,l)$ for all $l < k$, the 
conditional probability of $E(H,v,w,k)$ given $\F_{t+k-1}$
can be bounded below: inductively, Lemma~\ref{lem:min}
bounds from below the product of $W(v_j , v_{j-1} , t) W(v_j , v_{j+k} , t)$,
and hence the probability that $U(v_j , t) = \{ v_j , v_{j-1} , v_{j+k} \}$;
these conditionally independent probabilities may then be multiplied
to prove the claim, with the bound depending only on $x$ and $N$.

>From this argument, we see that the intersection $E(H,v,w) := 
\bigcap_{1 \leq k \leq r-2} E(H,v,w,k)$ has a probability which 
is bounded from below.  Sequentially,
we may choose a sequence of values for $w$ running through all 
vertices of $H$ at some distance $r(w) - 1 \geq 2$ from $v$, measured 
in the metric on $H$.  For each such $w$, we can bound from below 
the probability that in $r-2$ more time steps the path from $v$ to $w$ 
will be transitively completed.  We denote these events $E'(H,v,w)$, the
prime denoting the time shift to allow events analogous to $E(H,v,w)$
to occur sequentially.  Summing the time to run over all $w \in H$
yields at most $N^2$ time steps.   Let $E(H,v)$ denote the intersection 
of all the events $E' (H,v,w)$.   Inductively, we see that the probability
of $E(H,v)$ is bounded from below by a positive number depending only on 
$N$ and $x$.  

Finally, we let $(H,v)$ vary with $H$ exhausting components
of $G(t)$ and $v$ a choice function on the vertices of $H$.   
The events $E(H,v)$ are all conditionally independent given $\F_t$,
so the probability of their intersection, $E$, is bounded from below by 
a positive constant which we call $c$.  By Lemma~\ref{lem:min} once more,
on $E$, we know that~(\ref{eq:min}) is satisfied for each $e \in 
\overline{G(t)}$.    $\Cox$

\noindent{\sc Proof of Theorem}~\ref{th:limit}:
For any subset $V$ of agents, let 
$$E(V,t) := \bigcap_{s \geq t} \bigcap_{v \in V, w \in V^c} \left \{
    \{ v , w \} \notin G(s) \right \}$$
denote the event that from time $t$ onward, $V$ is isolated from 
its complement.  If $V$ is the vertex set of a component of
$G(t)$, then the conditional probability given $\F_t$ of the 
event $E(V,t)$ may be bounded from below as follows.  For any
$v \in V , w \in V^c$, and for any $s \geq t$, if the edge 
$e := \{ v , w \}$ is not in $G(r)$ for any $t \leq r < s$, then
by part~1 of Lemma~\ref{lem:min}, its weight $W(e,s)$ is at most 
$(1-x)^{s-t} 3 N x^{-1}$.  Since $\sum_z W(v,z,s) \geq 2$ for
all $v,z,s$, it follows from the evolution equations that
$$\P (e \in G(s) \| F_s) \leq \frac{(1-x)^{s-t} 3 N x^{-1}} 
   {2 + (1-x)^{s-t} 3 N x^{-1}} \, .$$
It follows that 
$$\P \left ( \exists v \in V , w \in V^c \, : \{ v , w \} \in G(s) \| \F_s
   \right ) = O \left ( N x^{-1} (1-x)^{s-t} \right ) $$
uniformly in $N , x$ and $t$ as $s - t \to \infty$ (though the uniformity
in $N$ and $x$ is not needed).  By the conditional Borel-Cantelli Lemma,
it follows that
\begin{equation} \label{eq:V}
\P (E(V,t) \| F_t) > c(N,x)
\end{equation}
on the event that $V$ is the vertex set of a component 
of $\overline{G(t)}$.  

By the reverse direction of the Conditional Borel-Cantelli Lemma, 
the event $E(V,t)$ occurs for some $t$ with probability~1
on that event that $V$ is a component of $\overline{G(t)}$ 
infinitely often.  Let $e = \{ v , w \}$ be any edge.  If $e \notin 
\overline{G(t)}$ infinitely often, then since there are only finitely 
many subsets of vertices, it follows that $v \in V$ and $w \in W$
for some disjoint $V$ and $W$ that are infinitely often components
of $G(t)$.  This implies that $e \in \overline{G(t)}$ finitely often.
We have shown that, almost surely, the edges come in two types: those in 
$\overline{G(t)}$ finitely often and those in $\overline{G(t)}$
all but finitely often.  This further implies that $\overline{G(t)}$
is eventually constant.  Denote this almost sure limit by $G_\infty$.
It remains to characterize $G_\infty$.

It is evident that $G_\infty$ contains no component of size less than
three, since $G(t)$ is the union of triangles $U(i,t)$.  
Suppose that $\overline{G(t)} = H$ for some $H$ of cardinality at least 
six.  By Lemma~\ref{lem:completion}, conditional on $\F_t$ 
and $\overline{G(t)} = H$, 
$$W(e,t+N^2) \geq \frac{(1-x)^{N^2}}{3N}$$
for every $e \in H$.  Write $H$ as the disjoint union of sets $J$ and $K$,
each of cardinality at least three.  Then with probability at least
$$\left ( \frac{(1-x)^{N^2}}{3N + (1-x)^{N^2}} \right )^{|J| + |K|}$$
$U(i , t + N^2) \subseteq J$ for every $i \in J$ and $U(i , t + N^2) 
\subseteq K$ for every $i \in K$.  In this case, $\overline{G(t+N^2)}$
has components that are proper subsets of $H$.  By the martingale
convergence theorem, 
$$\P (H \mbox{ is a component of } G_\infty \| \F_t)$$
converges with probability~1 to the indicator function of $H$ being 
a component of $G_\infty$.  From the above computation, it is not 
possible for $\P (H \mbox{ is a component of } G_\infty \| \F_t)$ 
to converge to~1 when $H$ has cardinality six or more.  Therefore, 
every component of $G_\infty$ has cardinality~3, 4 or~5.

The rest of the proof is easy.  Let $V_1 , \ldots , V_k$ be any
partition of $[N]$ into sets of cardinalities 3, 4 and~5.  The
derivation of~(\ref{eq:V}) shows that
$$\P \left ( \bigcap_{j=1}^k E(V_j , 1) \right ) > 0 \, ,$$
in other words, with positive probability $G_\infty$ has $k$
components which are precisely the complete graphs on $V_1 , \ldots , V_k$.  
It is elementary that a coupling may be produced between the Three's 
Company processes on populations of sizes $N$ and $K < N$ (with the 
same $x$ value), so that if $\{ \Wt (i,j,t) , \Ut (i,t) \}$ 
are the weight and choice variables for the smaller population, 
then $\Ut (i,t) = U(i,t)$ and $\Wt (i,j,t+1) = W(i,j,t+1)$ for all 
$t < \tau$ where $\tau$ is the first time, possibly infinite, at which 
$U(i,t)$ contains an edge between $[K]$ and $\{ K+1 , \ldots , N \}$.
In general, coupling methods show that if $\P (G_\infty = G_0 \| t) > 1 - \ee$
then the conditional distribution of the Three's Company process from
time $t$ onward given $\F_t$ and $G_\infty = G_0$, shifted back 
$t$ time units and restricted to a component $H$ of $G_\infty$, 
is within $\ee$ in total variation of the distribution of 
the Three's Company process on $H$ started with initial weights
$W' (i,j,0) := W(i,j,t)$.  

The Three's company process on a population of size 3, 4 or 5 
with any discount rate $1-x < 1$ is ergodic: to see this just note
that the Markov chain whose state space is the collection of $W$
variables is Harris recurrent as a consequence of 
Lemma~\ref{lem:completion}.  The invariant measure gives positive
weight to each edge, so each agent chooses each other with positive
frequency, finishing the proof of Theorem~\ref{th:limit}.   $\Cox$

\end{document}